\documentclass[global]{svjour}
\usepackage{amssymb}
\usepackage{amsfonts}
\usepackage{amsmath}

\input{tcilatex}

\begin{document}

\title{On the existence of solutions for a generalized strong vector
quasi-equilibrium problem}
\author{Monica Patriche}
\institute{University of Bucharest 
\email{monica.patriche@yahoo.com}%
}
\mail{\\
University of Bucharest, Faculty of Mathematics and Computer Science, 14
Academiei Street, 010014 Bucharest, Romania}
\maketitle

\bigskip \textbf{Abstract: }In this paper, we consider a generalized strong
vector quasi-equilibrium problem and we prove the existence of its solutions
by using some suxiliary results. One of the established theorems is proved
by using an approximation method.



\section{Introduction}

Vector equilibrium problem is a unified model of several problems, for
instance, vector variational inequalities and vector optimization problems.
For further relevant information on this topic, the reader is referred to
the following recent publications available in our bibliography: [1-7],
[9-12], [14-17].

In this paper, we will suppose that $X$\textit{\ }is a nonempty, convex and
compact set in a Hausdorff locally convex space $E$ , $A:X\rightarrow 2^{X}$
and $F:X\times X\times X\rightarrow 2^{X}$ are correspondences and $C\subset
X$ is a nonempty closed convex cone with int$C\neq \emptyset $.

We consider the following generalized strong vector quasi-equilibrium
problem (in short, GSVQEOP):\newline

find $x^{\ast }\in X$ such that $x^{\ast }\in \overline{A}(x^{\ast })$ and
each $u\in A(x^{\ast })$ implies that $F(u,x^{\ast },z)\nsubseteq $int$C$
for each $z\in A(x^{\ast }).$

\section{Preliminary results}

Let $X$, $Y$ be topological spaces and $T:X\rightarrow 2^{Y}$ be a
correspondence. $T$ is said to be \textit{upper semicontinuous} if for each $%
x\in X$ and each open set $V$ in $Y$ with $T(x)\subset V$, there exists an
open neighborhood $U$ of $x$ in $X$ such that $T(x)\subset V$ for each $y\in
U$. $T$ is said to be \textit{lower semicontinuous} if for each x$\in X$ and
each open set $V$ in $Y$ with $T(x)\cap V\neq \emptyset $, there exists an
open neighborhood $U$ of $x$ in $X$ such that $T(y)\cap V\neq \emptyset $
for each $y\in U$. $T$ is said to have \textit{open lower sections} if $%
T^{-1}(y):=\{x\in X:y\in T(x)\}$ is open in $X$ for each $y\in Y.$

The following lemma will be crucial in the proofs.

\begin{lemma}
(Yannelis and Prabhakar, \cite{yan}). \textit{Let }$X$\textit{\ be a
paracompact Hausdorff topological space and }$Y$\textit{\ be a Hausdorff
topological vector space. Let }$T:X\rightarrow 2^{Y}$\textit{\ be a
correspondence with nonempty convex values\ and for each }$y\in Y$\textit{, }%
$T^{-1}(y)$\textit{\ is open in }$X$\textit{. Then, }$T$\textit{\ has a
continuous selection that is, there exists a continuous function }$%
f:X\rightarrow Y$\textit{\ such that }$f(x)\in T(x)$\textit{\ for each }$%
x\in X$\textit{.\medskip }
\end{lemma}

The correspondence $\overline{T}$ is defined by $\overline{T}(x):=\{y\in
Y:(x,y)\in $cl$_{X\times Y}$ Gr $T\}$ (the set cl$_{X\times Y}$ Gr $(T)$ is
called the adherence of the graph of $T$)$.$ It is easy to see that cl $%
T(x)\subset \overline{T}(x)$ for each $x\in X.$

If $X$ and $Y$ are topological vector spaces, $K$ is a nonempty subset of $%
X, $ $C$ is a nonempty closed convex cone and $T:K\rightarrow 2^{Y}$ is a
correspondence, then \cite{luc}, $T$ is called \textit{upper }$C$\textit{%
-continuous at} $x_{0}\in K$ if, for any neighborhood $U$ of the origin in $%
Y,$ there is a neighborhood $V$ of $x_{0}$ such that, for all $x\in V,$ $%
T(x)\subset T(x_{0})+U+C.$ $T$ is called \textit{lower }$C$\textit{%
-continuous at} $x_{0}\in K$ if, for any neighborhood $U$ of the origin in $%
Y,$ there is a neighborhood $V$ of $x_{0}$ such that, for all $x\in V,$ $%
T(x_{0})\subset T(x)+U-C.$

The property of properly $C-$quasi-convexity for correspondences is
presented below.

Let $X$ be a nonempty convex subset of a topological vector space\textit{\ }$%
E,$ $Y$ be a topological vector space, and $C$ be a pointed closed convex
cone in $Z$ with its interior int$C\neq \emptyset .$ Let $T:X\rightarrow
2^{Y}$ be a correspondence with nonempty values. $T$ is said to be \textit{%
properly }$C-$\textit{quasi-convex on} $X$ (\cite{long}), if for any $%
x_{1},x_{2}\in X$ and $\lambda \in \lbrack 0,1],$ either $T(x_{1})\subset
T(\lambda x_{1}+(1-\lambda )x_{2})+C$ or $T(x_{2})\subset T(\lambda
x_{1}+(1-\lambda )x_{2})+C.\medskip $

In order to establish our main theorems, we need to prove some auxiliary
results. The starting point is the following statement:

\begin{theorem}
Let $X$ be \textit{a nonempty, convex and compact set in a} \textit{%
Hausdorff locally convex space }$E$ and let $\mathcal{C\ }$\ be a nonempty
subset of $X\times X.$\textit{\ Assume that t}he following conditions are
fulfilled:
\end{theorem}

\textit{a) }$\mathcal{C}^{-}(y)$ $=\{x\in X:(x,y)\in C\}$ \textit{is open for%
} \textit{each } $y\in X;$

\textit{b)} $\mathcal{C}^{+}(x)=\{y\in X:(x,y)\in C\}\mathit{\ }$ \textit{is
convex and nonempty for each} $x\in X.$

\textit{Then, there exists }$x^{\ast }\in X$ such that $(x^{\ast },x^{\ast
})\in C$\textit{.\medskip }

\begin{proof}
Let us define the correspondence $T:X\rightarrow 2^{X},$ by

$T(x)=\mathcal{C}^{+}(x)$ for each $x\in X.$

The correspondence $T$ is nonempty and convex valued and it has open lower
sections$.$

We apply the Yannelis and Prabhakar's Lemma and we obtain that $T$\textit{\ }%
has a continuous selection $f:X\rightarrow X.$

According to the Tychonoff fixed point Theorem \cite{is}, there exists $%
x^{\ast }\in X$ such that $f(x^{\ast })=x^{\ast }.$ Hence, $x^{\ast }\in
T(x^{\ast })$ and obviously, $(x^{\ast },x^{\ast })\in \mathcal{C}.\medskip
\medskip $
\end{proof}

The next two results are direct consequences of Theorem 1.\medskip

\begin{theorem}
Let $X$ be \textit{a nonempty, convex and compact set in a} \textit{%
Hausdorff locally convex space }$E$ , and let $A:X\rightarrow 2^{X}$ and $%
P:X\times X\rightarrow 2^{X}$ be correspondences such that the following
conditions are fulfilled:
\end{theorem}

\textit{a)} $A$\textit{\ has nonempty, convex values and open lower sections;%
}

\textit{b)} \textit{the set }$\{y\in X:$ $A(x)\cap P(x,y)=\emptyset \}\cap
A(x)$ \ \textit{is nonempty for each} $x\in X;$

\textit{c) }$\{y\in X:$\textit{\ }$A(x)\cap P(x,y)=\emptyset \}$\textit{\ is
convex for each }$x\in X;$

\textit{d)} $\{x\in X:A(x)\cap P(x,y)=\emptyset \}\mathit{\ }$ \textit{is
open for each} $y\in X.$

\textit{Then, there exists }$x^{\ast }\in $\textit{\ }$X$\textit{\ such that 
}$\ x^{\ast }\in A(x^{\ast })$ \textit{and }$A(x^{\ast })\cap P(x^{\ast
},x^{\ast })=\emptyset .\medskip $

\begin{proof}
Let us define the set $\mathcal{C}=\{(x,y)\in X\times X:$ $A(x)\cap
P(x,y)=\emptyset \}\cap $Gr$A.$

Then, $\mathcal{C}^{+}(x)=\{y\in X:$ $A(x)\cap P(x,y)=\emptyset \}\cap A(x)$
for each $\ x\in X$ and

$\mathcal{C}^{-}(y)=A^{-1}(y)\cap \{x\in X:A(x)\cap P(x,y)=\emptyset \}\ $%
for each $x\in X.$

Assumption b) implies that $\mathcal{C}$ is nonempty$.$ The set $\mathcal{C}%
^{-}(y)$ is open for each $y\in X$ since Assumptions a) and d) hold.

According to Assumptions b) and c), $A(x)\cap $ $\mathcal{C}^{+}(x)$ is
nonempty and convex for each $x\in X.$

All hypotheses of Theorem 1 are fulfilled, and then, there exists $x^{\ast
}\in $\ $X$\ such that $\ x^{\ast }\in A(x^{\ast })$ and $A(x^{\ast })\cap
P(x^{\ast },x^{\ast })=\emptyset .$
\end{proof}

We establish the following result as a consequence of Theorem 2. It will be
used in order to prove the existence of solutions for the considered vector
quasi-equilibrium problem.

\begin{theorem}
Let $X$ be \textit{a nonempty, convex and compact set in a} \textit{%
Hausdorff locally convex space }$E$ , and let $A:X\rightarrow 2^{X}$, $%
P:X\times X\rightarrow 2^{X}$ be correspondences such that the following
conditions are fulfilled:
\end{theorem}

\textit{a)} $A$\textit{\ has nonempty, convex values and open lower sections;%
}

\textit{b)} \textit{the set }$\{y\in X:$ $\overline{A}(x)\cap
P(x,y)=\emptyset \}\cap A(x)$ \ \textit{is nonempty for each} $x\in X;$

\textit{c) }$\{y\in X:$\textit{\ }$\overline{A}(x)\cap P(x,y)=\emptyset \}$%
\textit{\ is convex for each }$x\in X;$

\textit{d)} $\{x\in X:\overline{A}(x)\cap P(x,y)=\emptyset \}\mathit{\ }$ 
\textit{is open for each} $y\in X.$

\textit{Then, there exists }$x^{\ast }\in $\textit{\ }$X$\ \textit{\ such
that }$\ x^{\ast }\in A(x^{\ast })$ \textit{and }$A(x^{\ast })\cap P(x^{\ast
},x^{\ast })=\emptyset .\medskip $

We note that, according to Theorem 2, there exists $x^{\ast }\in $\textit{\ }%
$X$\textit{\ }such that $\ x^{\ast }\in A(x^{\ast })$ and\textit{\ }$%
\overline{A}(x^{\ast })\cap P(x^{\ast },x^{\ast })=\emptyset .$ Obviously, $%
\overline{A}(x^{\ast })\cap P(x^{\ast },x^{\ast })=\emptyset $ implies $%
A(x^{\ast })\cap P(x^{\ast },x^{\ast })=\emptyset .\medskip $

If $A(x)=X$ for each $x\in X$, Theorem 2 implies the following corollary.

\begin{corollary}
Let $X$ be \textit{a nonempty, convex and compact set in a} \textit{%
Hausdorff locally convex space }$E$ , and let $P:X\times X\rightarrow 2^{X}$
be a correspondence such that the following conditions are fulfilled:
\end{corollary}

\textit{a) }$\{y\in X:$\textit{\ }$P(x,y)=\emptyset \}$\textit{\ is nonempty
and convex for each }$x\in X;$

\textit{b)} $\{x\in X:P(x,y)=\emptyset \}\mathit{\ }$ \textit{is open for
each} $y\in X.$

\textit{Then, there exists }$x^{\ast }\in $\textit{\ }$X$\textit{\ such that 
}$\ P(x^{\ast },x^{\ast })=\emptyset .\medskip $

By applying an approximation method of proof, we can prove Theorem 4.

\begin{theorem}
Let $X$ be \textit{a nonempty, convex and compact set in a} \textit{%
Hausdorff locally convex space }$E$ and $\mathcal{C\ }$\ be a nonempty,
closed subset of $X\times X.$\textit{\ Assume that there exists a sequence }$%
(G_{k})_{k\in \mathbb{N}^{\ast }}$ \textit{of subsets of }$X\times X$ such
that \textit{t}he following conditions are fulfilled:
\end{theorem}

\textit{a) }for each $k\in \mathbb{N}^{\ast },$ $G_{k}^{-}(y)$ $=\{x\in
X:(x,y)\in G_{k}\}$ \textit{is open for} \textit{each } $y\in X;$

\textit{b)} for each $k\in \mathbb{N}^{\ast },$ $G_{k}^{+}(x)=\{y\in
X:(x,y)\in G_{k}\}\mathit{\ }$ \textit{is convex and nonempty for each} $%
x\in X;$

\textit{c) }$G_{k}\supseteq G_{k+1}$\textit{\ for each }$k\in \mathbb{N}%
^{\ast };$

\textit{d) for every open set }$G$\textit{\ with }$G\supset \mathcal{C},$%
\textit{\ there exists }$k\in \mathbb{N}^{\ast }$\textit{\ such that }$%
G_{k}\subseteq \mathcal{C}.$

\textit{Then, there exists }$x^{\ast }\in X$\textit{\ such that }$(x^{\ast
},x^{\ast })\in \mathcal{C}$\textit{.\medskip }

\begin{proof}
For each $k\in \mathbb{N}^{\ast },$ we apply Theorem 1. Let $x^{k}\in X$
such that $(x^{k},x^{k})\in G_{k}.$ Since $X$ is a compact set, we can
consider that the sequence $(x^{k})_{k}$ converges to some $x^{\ast }\in X.$
We claim that $(x^{\ast },x^{\ast })\in \mathcal{C}.$

Indeed, let us suppose, by way of contradiction, that $(x^{\ast },x^{\ast
})\notin \mathcal{C}.$ Since $\mathcal{C}$\textit{\ }is nonempty and
compact, we can choose a neighborhood $V_{(x^{\ast },x^{\ast })}$ of $%
(x^{\ast },x^{\ast })$ and an open set $G$ such that $G\supset \mathcal{C}$
and $V_{(x^{\ast },x^{\ast })}\cap G=\emptyset .$ According to Assumptions
d) and c), there exists $k_{1}\in \mathbb{N}^{\ast }$ such that $%
G_{k}\subseteq G$ for each $k\geq k_{1}.$ Since $V_{(x^{\ast },x^{\ast })}$
is a neighborhood of $(x^{\ast },x^{\ast }),$ there exists $k_{2}\in \mathbb{%
N}^{\ast }$ such that $(x^{k},x^{k})\in V_{(x^{\ast },x^{\ast })}$ for each $%
k\geq k_{2}.$ $\ $Hence, for $k\geq $max($k_{2},k_{1}),$ $%
(x^{k},x^{k})\notin G_{k},$ which is a contradiction.

Consequently, $(x^{\ast },x^{\ast })\in \mathcal{C}.$
\end{proof}

Theorem 5 is a consequence of Theorem 4 and it will be used in Section 3 in
order to prove the existence of solutions for GSVQEP.

\begin{theorem}
\textit{Let }$X$\textit{\ be a nonempty, convex and compact set in a
Hausdorff locally convex space }$E$\textit{\ , and let }$A:X\rightarrow
2^{X} $\textit{\ and }$P:X\times X\rightarrow 2^{X}$\textit{\ be
correspondences such that the following conditions are fulfilled:}
\end{theorem}

\textit{a) }$A$\textit{\ has nonempty, convex values and open lower sections;%
}

\textit{b) the set }$U=\{(x,y)\in X\times X:$\textit{\ }$\overline{A}(x)\cap
P(x,y)=\emptyset \}$\textit{\ \ is closed} \textit{and} $U\cap $Gr$\overline{%
A}$ \textit{is nonempty;}

\textit{c) there exists a sequence }$(P_{k})_{k\in N}$\textit{\ of
correspondences, where, for each }$k\in \mathbb{N}^{\ast },$\textit{\ }$%
P_{k}:X\times X\rightarrow 2^{X}$\textit{\ and let }$U_{k}=\{(x,y)\in
X\times X:$\textit{\ }$A(x)\cap P_{k}(x,y)=\emptyset \}$\textit{. Assume
that:}

\textit{\ }$\ \ \ $\textit{c1) }$U_{k}^{+}(x)=\{y\in X:$\textit{\ }$%
\overline{A}(x)\cap P_{k}(x,y)=\emptyset \}$\textit{\ is convex for each }$%
x\in X$\textit{\ and }$U_{k}^{+}(x)\cap A(x)\neq \emptyset $\textit{\ for
each }$x\in X;$

\textit{\ \ \ \ c2 ) }$U_{k}^{-}(y)=\{x\in X:\overline{A}(x)\cap
P_{k}(x,y)=\emptyset \}\ $\textit{\ is open for each }$y\in X;$

\textit{\ \ \ \ c3) for each }$k\in \mathbb{N}^{\ast },$\textit{\ }$%
P_{k}(x,y)\subseteq P_{k+1}(x,y)$\textit{\ for each }$(x,y)\in X\times X;$

\textit{\ \ \ \ \ c4) for every open set }$G$\textit{\ with }$G\supset U\cap 
$Gr$\overline{A},$\textit{\ there exists }$k\in \mathbb{N}^{\ast }$\textit{\
such that }$G\supseteq U_{k}\cap $Gr$A.$

\textit{Then, there exists }$x^{\ast }\in $\textit{\ }$X$\textit{\ such that 
}$\ x^{\ast }\in \overline{A}(x^{\ast })$\textit{\ and }$A(x^{\ast })\cap
P(x^{\ast },x^{\ast })=\emptyset .\medskip $

\begin{proof}
Let us define $\mathcal{C=}U\cap $Gr$\overline{A}.$ According to Assumptions
b) and c), $\mathcal{C\ }$\ is a nonempty and closed subset of $X\times X.$

Further, \textit{\ }for each\textit{\ }$k\in \mathbb{N}^{\ast },$ let us
define $G_{k}=U_{k}\cap $Gr$A\subseteq X\times X.$

Then, for each $k\in \mathbb{N}^{\ast },$ $G_{k}^{+}(x)$ $=\{y\in X:(x,y)\in
G_{k}\}=U_{k}^{+}(x)\cap A(x)$ is nonempty and convex for each\textit{\ }$%
x\in X,$ since Assumptions a) and c1) hold.

For each $k\in \mathbb{N}^{\ast },$ $G_{k}^{-}(y)$ $=\{x\in X:(x,y)\in
G_{k}\}=U_{k}^{-}(y)\cap A^{-1}(y)$ is open for each\textit{\ }$y\in X,$
since Assumptions a) and c2) hold.

Assumption c3) implies that, \textit{\ }for each\textit{\ }$k\in \mathbb{N}%
^{\ast },$ $U_{k+1}\subseteq U_{k}$ and then, $G_{k}\supseteq G_{k+1}$ and
Assumption c4) implies that for every open set $G$ with $G\supset \mathcal{C}%
,$ there exists $k\in \mathbb{N}^{\ast }$ such that $G_{k}\subseteq \mathcal{%
C}.$

All hypotheses of Theorem 4 are verified. Therefore, there exists $x^{\ast
}\in X$ such that $(x^{\ast },x^{\ast })\in \mathcal{C}$.

Consequently, there exists $x^{\ast }\in X$ such that $\ x^{\ast }\in 
\overline{A}(x^{\ast })$\textit{\ }and\textit{\ }$A(x^{\ast })\cap P(x^{\ast
},x^{\ast })=\emptyset .\medskip $
\end{proof}

\section{Main results}

This section is devoted to the study of the existence of solutions for the
considered generalized strong vector quasi-equilibrium problem. We derive
our main results by using the auxiliary theorems concerning correspondences,
which have been established in the previous section. This new approach to
solve GSVQEP is intended to provide new conditions under which the solutions
exist.\medskip

The first theorem states that GSVQEP has solutions if $F(\cdot ,y,\cdot )$
is lower ($-C$)-semicontinuous for each $y\in X$ and $F(u,\cdot ,z)$\ is
properly $C-$\ quasi-convex for each $(u,z)\in X\times X.$

\begin{theorem}
\textit{Let }$F:X\times X\times X\rightarrow 2^{X}$\textit{\ be a
correspondence with nonempty values. Suppose that:}
\end{theorem}

\textit{a) }$A$\textit{\ has nonempty, convex values and open lower sections;%
}

\textit{b) for each }$x\in X,$\textit{\ there exists }$y\in A(x)$\textit{\
such that each }$u\in \overline{A}(x)$\textit{\ implies that }$%
F(u,y,z)\nsubseteq $\textit{int}$C$\textit{\ for each }$z\in \overline{A}(x);
$

\textit{c) }$F(\cdot ,y,\cdot ):$\textit{\ }$X\times X\rightarrow 2^{X}$%
\textit{\ is lower (}$-C$\textit{)-semicontinuous for each }$y\in X;$

\textit{d) for each }$(u,z)\in X\times X,$\textit{\ }$F(u,\cdot
,z):X\rightarrow 2^{X}$\textit{\ is properly }$C-$\textit{\ quasi-convex.}

\textit{Then, there exists }$x^{\ast }\in X$\textit{\ such that }$x^{\ast
}\in A(x^{\ast })$\textit{\ and each }$u\in A(x^{\ast })$\textit{\ implies
that }$F(u,x^{\ast },z)\nsubseteq $\textit{int}$C$\textit{\ for each }$z\in
A(x^{\ast }),$\textit{\ that is, }$x^{\ast }$\textit{\ is a solution for
GSVQEP.\medskip }

\begin{proof}
Let us define $P:X\times X\rightarrow 2^{X},$ by

$P(x,y)=\{u\in X:$ $\exists z\in \overline{A}(x)$ such that $%
F(u,y,z)\subseteq C\}$ for each $(x,y)\in X\times X.$

Assumption b) implies that the set $\{y\in X:$ $\overline{A}(x)\cap
P(x,y)=\emptyset \}\cap A(x)$ \ is nonempty for each $x\in X.$

We claim that the set $E(x)$ $\ $is convex for each $x\in X,$ where\ $%
E(x)=\{y\in X:$\textit{\ }$\overline{A}(x)\cap P(x,y)=\emptyset \}=$

$=\{y\in X:$\textit{\ }each $u\in \overline{A}(x)$ implies that $%
F(u,y,z)\nsubseteq C$ for each $z\in \overline{A}(x)\}.$

Indeed, let us fix $x_{0}\in X$ and let us consider $y_{1},y_{2}\in
E(x_{0}). $ This means that each $u\in \overline{A}(x_{0})$ implies that $%
F(u,y_{1},z)\nsubseteq C$ and $F(u,y_{2},z)\nsubseteq C$ for each $z\in 
\overline{A}(x_{0}).$

Let $y(\lambda )=\lambda y_{1}+(1-\lambda )y_{2}$ be defined for each $%
\lambda \in \lbrack 0,1].$

We claim that $y(\lambda )\in E(x_{0})$ for each $\lambda \in \lbrack 0,1].$

Suppose, on the contrary, that there exist $\lambda _{0}\in \lbrack 0,1]$ $,$
$u^{\prime }\in \overline{A}(x_{0})$ and $z^{\prime }\in \overline{A}(x_{0})$
such that $F(u^{\prime },y(\lambda _{0}),z^{\prime })\subseteq C.$ Since $%
(F(u^{\prime },\cdot ,z^{\prime }):X\rightarrow 2^{X}$ is\textit{\ }properly%
\textit{\ }$C-$\textit{\ }quasi-convex$,$ we have that:

$F(u^{\prime },y_{1},z^{\prime })\subseteq F(u^{\prime },y(\lambda
),z^{\prime })+C$ or $F(u^{\prime },y_{2},z^{\prime })\subseteq F(u^{\prime
},y(\lambda ),z^{\prime })+C$

On the other hand, it is true that $F(u^{\prime },y(\lambda ),z^{\prime
})\subseteq C.$ We obtain that:

$F(u^{\prime },y_{j},z^{\prime })\subseteq C+C\subseteq C$ for $j=1$ or for $%
j=2$.

This contradicts the assumption that $y_{1},y_{2}\in E(x_{0})$.
Consequently, $E(x_{0})$ is convex and Assumption c) from Theorem 3 is
fulfilled.

Now, we will prove that $D(y)=\{x\in X:\overline{A}(x)\cap P(x,y)=\emptyset
\}\mathit{\ }$ is open for each $y\in X.$

In order to do this, we will show that $^{C}D(y)$ is closed for each $y\in X,
$ where $\ ^{C}D(y)=\{x\in X:\overline{A}(x)\cap P(x,y)\neq \emptyset \}%
\mathit{=}\{x\in X:$ there exist $u,z\in \overline{A}(x)$ such that $%
F(u,y,z)\subseteq C\}.$

Let $(x_{\alpha })_{\alpha \in \Lambda }$ be a net in $^{C}D(y)$ such that
lim$_{\alpha }x_{\alpha }=x_{0}.$ Then, there exist $u_{\alpha },z_{\alpha
}\in \overline{A}(x_{\alpha })$ such that $F(u_{\alpha },y,z_{\alpha
})\subseteq C.$

Since $X$ is a compact set, we can suppose that $(u_{\alpha })_{\alpha \in
\Lambda },(z_{\alpha })_{\alpha \in \Lambda }$ are convergent nets and let
lim$_{\alpha }u_{\alpha }=u_{0}$ and lim$_{\alpha }z_{\alpha }=z_{0}.$

The closedness of $\overline{A}$ implies that $u_{0},z_{0}\in \overline{A}%
(x_{0}).$

Now, we claim that $F(u_{0},y,z_{0})\subseteq C.$

Since $F(u_{\alpha },y,z_{\alpha })\subseteq C$ and $F(\cdot ,y,\cdot ):$%
\textit{\ }$X\times X\rightarrow 2^{X}$\textit{\ }is lower ($-C$%
)-semicontinuous, for each neighborhood $U$ of the origin in $X,$ there
exists a subnet $(u_{\beta },z_{\beta })_{\beta }$ of $(u_{\alpha
},z_{\alpha })_{\alpha }$ such that $F(u_{0},y,z_{0})\subset F(u_{\beta
},y,z_{\beta })+U+C.$ Hence, $F(u_{0},y,z_{0})\subset U+C.$

We will show that $F(u_{0},y,z_{0})\subset C.$ Suppose, by way of
contradiction, that there exists $t\in F(u_{0},y,z_{0})\cap ^{C}C.$ We note
that $B=C-t$ is a closed set which does not contain $0.$ It follows that $%
^{C}B$ is open and contains $0.$ Since $X$ is locally convex, there exists a
convex neighborhood $U_{1}$ of origin such that $U_{1}\subset X\backslash B$
and $U_{1}=-U_{1}$. Thus, $0\notin B+U_{1}$ and then, $t\notin C+U_{1},$
which is a contradiction. Therefore, $F(u_{0},y,z_{0})\subset C.$

We proved that there exist $u_{0},z_{0}\in \overline{A}(x_{0})$ such that $%
F(u_{0},y,z_{0})\subseteq C.$ It follows that $^{C}D(y)$ is closed. Then, $%
D(y)$ is an open set and Assumption d) from Theorem 3 is fulfilled.

Consequently, all conditions of Theorem 3 are verified, so that there exists 
$x^{\ast }\in $\ $X$\ such that $\ x^{\ast }\in A(x^{\ast })$ and $A(x^{\ast
})\cap P(x^{\ast },x^{\ast })=\emptyset .\medskip $ Obviously, $x^{\ast }$%
\textit{\ }is a solution for GSVQEP.\medskip
\end{proof}

In order to obtain a second result concerning the existence of solutions of
GSVQEP, we use an approximation method and Theorem 5. We mention that this
result does not require convexity properties for the correspondence $F.$

\begin{theorem}
\textit{Let }$F:X\times X\times X\rightarrow 2^{X}$\textit{\ be a
correspondence. Suppose that:}
\end{theorem}

\textit{a) }$A$\textit{\ has nonempty, convex values and open lower
sections; }$\overline{A}$\textit{\ is lower semicontinuous;}

\textit{b) }$F$\textit{\ is upper semicontinuous with nonempty, closed
values;}

\textit{c) }$U\cap $Gr$\overline{A}$\textit{\ is nonempty, where }$U=\{$%
\textit{\ }$(x,y)\in X\times X:u\in \overline{A}(x)$\textit{\ implies that }$%
F(u,y,z)\nsubseteq $\textit{int}$C$\textit{\ for each }$z\in \overline{A}%
(x)\};$

\textit{d) there exists a sequence }$(F_{k})_{k\in N}$\textit{\ of
correspondences, such that, for each }$k\in \mathbb{N}^{\ast },$\textit{\ }$%
F_{k}:X\times X\times X\rightarrow 2^{X}$\textit{\ and let }$%
U_{k}=\{(x,y)\in X\times X:u\in \overline{A}(x)$\textit{\ implies that }$%
F_{k}(u,y,z)\nsubseteq $\textit{int}$C$\textit{\ for each }$z\in \overline{A}%
(x)\}$\textit{. Assume that:}

\textit{d1) for each }$k\in \mathbb{N}^{\ast }$\textit{\ and for each }$x\in
X,$\textit{\ there exists }$y\in A(x)$\textit{\ such that each }$u\in 
\overline{A}(x)$\textit{\ implies that }$F_{k}(u,y,z)\nsubseteq $\textit{int}%
$C$\textit{\ for each }$z\in A(x);$

\textit{d2) for each }$k\in \mathbb{N}^{\ast }$\textit{\ and for each }$%
(u,z)\in X\times X,$\textit{\ }$F_{k}(u,\cdot ,z):X\rightarrow 2^{X}$\textit{%
\ is properly }$C-$\textit{\ quasi-convex;}

\textit{d3) for each }$k\in \mathbb{N}^{\ast }$ \textit{and} \textit{\ for
each }$y\in X,$\textit{\ }$F_{k}(\cdot ,y,\cdot ):$\textit{\ }$X\times
X\rightarrow 2^{X}$\textit{\ is lower (}$-C$\textit{)-semicontinuous}$;$

\textit{d4) for each }$k\in \mathbb{N}^{\ast },$\textit{\ for each }$%
(x,y)\in X\times X,$\textit{\ and for each }$u\in X$\textit{\ with the
property that }$\exists z\in \overline{A}(x)$\textit{\ such that }$%
F_{k}(u,y,z)\subseteq C,$\textit{\ there exists }$z^{\prime }\in \overline{A}%
(x)$\textit{\ such that }$F_{k+1}(u,y,z^{\prime })\subseteq C;$

\textit{d5) for every open set }$G$\textit{\ with }$G\supset U\cap $Gr$%
\overline{A},$\textit{\ there exists }$k\in \mathbb{N}^{\ast }$\textit{\
such that }$G\supseteq U_{k}\cap $\textit{Gr}$A.$

\textit{Then, there exists }$x^{\ast }\in X$\textit{\ such that }$x^{\ast
}\in \overline{A}(x^{\ast })$\textit{\ and each }$u\in A(x^{\ast })$\textit{%
\ implies that }$F(u,x^{\ast },z)\nsubseteq $\textit{int}$C$\textit{\ for
each }$z\in A(x^{\ast }),$\textit{\ that is, }$x^{\ast }$\textit{\ is a
solution for GSVQEP.\medskip }

\begin{proof}
Let us define $P:X\times X\rightarrow 2^{X},$ by

$P(x,y)=\{u\in X:$ $\exists z\in \overline{A}(x)$ such that $%
F(u,y,z)\subseteq C\}$ for each $(x,y)\in X\times X.$

We claim that, $U=\{$\textit{\ }$(x,y)\in X\times X:u\in \overline{A}(x)$%
\textit{\ }implies that $F(u,y,z)\nsubseteq $int$C$ for each $z\in \overline{%
A}(x)\}$ is closed$.$

Let $(x^{0},y^{0})\in $cl$U.$ Then, there exists $(x^{\alpha },y^{\alpha
})_{\alpha \in \Lambda }$ a net in $U$ such that $\lim_{\alpha }(x^{\alpha
},y^{\alpha })=(x^{0},y^{0})\in X\times X.$ Let $u\in \overline{A}(x^{0})$
and $z\in \overline{A}(x^{0}).$ Since $\overline{A}$ is lower semicontinuous
and $\lim_{\alpha }x^{\alpha }=x^{0},$ there exists nets $(u^{\alpha
})_{\alpha \in \Lambda }\ $\ and $(z^{\alpha })_{\alpha \in \Lambda }$ in $X$
such that $u^{\alpha },z^{\alpha }\in \overline{A}(x^{\alpha })$ for each $%
\alpha \in \Lambda $ and $\lim_{\alpha }u^{\alpha }=u,$ $\lim_{\alpha
}z^{\alpha }=z.$ Since $(x^{\alpha },y^{\alpha })_{\alpha \in \Lambda }$ is
a net in $U,$ $\ $then$,$ for each $\alpha \in \Lambda ,$ $F(u^{\alpha
},y^{\alpha },z^{\alpha })\nsubseteq $int$C$ $,$ that is, $F(u^{\alpha
},y^{\alpha },z^{\alpha })\cap W\neq \emptyset ,$ where $W=X\backslash $int$%
C,$ that is, there exists $(t^{\alpha })_{\alpha \in \Lambda }$ a net in $X$
such that $t^{\alpha }\in F(u^{\alpha },y^{\alpha },z^{\alpha })\cap W$ for
each $\alpha \in \Lambda .$

Since $X$ is compact, we can suppose that $\lim_{\alpha }t^{\alpha }=t^{0}.$
The closedness of $W$ implies that $t^{0}\in W.$ We invoke here the
closedness of $F$ and we conclude that $t^{0}\in F(u,y^{0},z).$ Therefore, $%
F(u,y^{0},z)\cap W\neq \emptyset ,$ and, thus, $u\in \overline{A}(x^{0})$
implies $F(u,y^{0},z)\nsubseteq $int$C$ for each $z\in \overline{A}(x^{0}).$
Hence, $U$ is closed.

For each $k\in \mathbb{N}^{\ast },$ let us define $P_{k}:X\times
X\rightarrow 2^{X},$ by

$P_{k}(x,y)=\{u\in X:$ $\exists z\in \overline{A}(x)$ such that $%
F_{k}(u,y,z)\subseteq C\}$ for each $(x,y)\in X\times X$ and

$U_{k}=\{(x,y)\in X\times X:u\in \overline{A}(x)$ implies that $%
F_{k}(u,y,z)\nsubseteq $int$C$ for each $z\in \overline{A}(x)\}=\{(x,y)\in
X\times X:$\textit{\ }$\overline{A}(x)\cap P_{k}(x,y)=\emptyset \}.$

\textit{\ }Let\textit{\ }$k\in \mathbb{N}^{\ast }.$

Assumption d1) implies that $U_{k}^{+}(x)\cap A(x)\neq \emptyset $\textit{\ }%
for each\textit{\ }$x\in X$ and Assumption d2) implies that $U_{k}^{+}(x)$
is convex \textit{\ }for each\textit{\ }$x\in X$ (we use a similar proof
with the one of Theorem 6)

Since $F_{k}(\cdot ,y,\cdot ):$\textit{\ }$X\times X\rightarrow 2^{X}$%
\textit{\ }is lower ($-C$)-semicontinuous for each $y\in X,$ by following an
argument similar with the one from the proof of Theorem 6, we can prove that

$U_{k}^{-}(y)=\{x\in X:\overline{A}(x)\cap P_{k}(x,y)=\emptyset \}\ $ is
open for each\textit{\ }$y\in X.$

Assumption d4) implies that $P_{k}(x,y)\subseteq P_{k+1}(x,y)$\textit{\ }for
each\textit{\ }$(x,y)\in X\times X.$

All conditions of Theorem 5 are verified, so that there exists $x^{\ast }\in 
$\ $X$\ such that $\ x^{\ast }\in \overline{A}(x^{\ast })$ and $A(x^{\ast
})\cap P(x^{\ast },x^{\ast })=\emptyset .\medskip $ Obvioulsy, $x^{\ast }$%
\textit{\ }is a solution for GSVQEP.\medskip
\end{proof}

\section{Concluding remarks}

This paper developed a framework for discussing the existence of solutions
for a generalized strong vector quasi-equilibrium problem. The results have
been obtained under assumptions which are different than the existing ones
in literature. An approximation technique of proof has been developed.

\end{document}